\documentclass[12pt,a4paper]{article}

\usepackage{amsmath}
\usepackage{amssymb}

\newtheorem{defn}{Definition}[section]
\newtheorem{thm}[defn]{Theorem}
\newtheorem{prop}[defn]{Proposition}
\newtheorem{lem}[defn]{Lemma}
\newtheorem{cor}[defn]{Corollary}

\newcommand\pf{\noindent{\bf Proof. }}
\newcommand\qed{\hfill\framebox(5,5){~}}
\newcommand\cf{\rm{c.f. }}

\newcommand{\g}{\mathfrak{g}}

\newcommand{\E}{{\cal{E}}}
\newcommand{\Flines}{{\cal{F}}}
\newcommand{\Pts}{{\cal{P}}}
\newcommand{\Lines}{{\cal{L}}}
\newcommand{\X}{{\cal{X}}}

\newcommand{\map}{\rightarrow}
\newcommand{\mapto}{\mapsto}

\newcommand{\dist}{{\rm{dist}}}

\newcommand{\rsh}{\rm{RSh}}

\newcommand{\Exp}{{\rm{Exp}}}
\newcommand{\aut}{{\rm{Aut}}}

\newcommand{\Z}{{\mathbb{Z}}}

\title{Recovering the Lie algebra from its extremal geometry}
\author{H. Cuypers, K. Roberts and S. Shpectorov}

\begin{document}
\maketitle

\begin{abstract} An element $x$ of a Lie algebra $L$ over the field $F$ is extremal if $[x,[x,L]]=Fx$. Under minor assumptions, it is known that, for a simple Lie algebra $L$, the extremal geometry $\E(L)$ is a subspace of the projective geometry of $L$ and either has no lines or is the root shadow space of an irreducible spherical building $\Delta$. We prove that if $\Delta$ is of simply-laced type, then $L$ is a quotient of a Chevalley algebra of the same type.
\end{abstract}

\section{Introduction}
Lie theory is an active and important field whose applications are widespread throughout mathematics and physics. Underlying objects in Lie theory are the simple Lie algebras which have striking connections to algebraic groups and buildings. Simple complex Lie algebras were classified by Killing \cite{k} and Cartan \cite{c} in the late nineteenth century. In the second half of the twentieth century, through the efforts of many mathematicians, a classification of all simple finite-dimensional Lie algebras over algebraically closed fields with characteristic at least $5$ was achieved (\cf \cite{s1,s2}). The two main classes of Lie algebras that emerge from the result are the classical Lie algebras derived from the simple complex Lie algebras and the Lie algebras of Cartan type derived from the Witt algebras. 

Given a simple complex Lie algebra $\g$ corresponding to a root system $\Phi$ of type $X_n$, we can construct its $\Z$-form $\g_\Z$ by taking the integral span of a Chevalley basis of $\g$. By tensoring $\g_\Z$ with an arbitrary field $F$ we obtain the Chevalley algebra $\g_F$. A \emph{classical Lie algebra} (or simply a \emph{classical algebra}) of type $X_n$ over $F$ is defined as an arbitrary nonzero quotient of the Chevalley algebra $\g_F$. For classical Lie algebras we have the familiar notions such as Cartan subalgebra, root elements, {\it etc}. 

Extremal elements are extensively utilised in the classification of simple modular Lie algebras over algebraically closed fields and there is good evidence to believe that they should play a crucial role for arbitrary fields.  An element $x$ of a Lie algebra $L$ over a field $F$ is \emph{extremal} if $[x,[x,L]] \subseteq Fx$ and such an element $x$ is a \emph{sandwich} if $[x,[x,L]] = 0$ (additional conditions are needed when $F$ has characteristic $2$). Examples of extremal elements are the long root elements of the classical algebras. In \cite{p}, Premet shows that any simple Lie algebra over an algebraically closed field of characteristic at least $5$ contains a nonzero extremal element. It is shown in \cite{cir} that, with one exception, a simple Lie algebra over an arbitrary field of characteristic at least $5$ that contains a non-sandwich extremal element is generated by such elements.

Extremal elements are also the focal point of the effort spearheaded by Cohen, \emph{et al.}, to provide a geometric characterisation of the classical Lie algebras. In \cite{ci1}, Cohen and Ivanyos define a point-line space on the set of nonzero extremal elements of $L$ called the \emph{extremal geometry} denoted by $\E(L)$. By combining the results in \cite{ci1, ci2, ks}, one can conclude that if $L$ is finite-dimensional, simple and generated by extremal elements with no sandwiches, then $\E(L)$ either contains no lines or is the root shadow space of a spherical building. An interesting and important question is whether one can recover the algebra $L$ from the geometry $\E(L)$. By using the theory of buildings first introduced by Tits in \cite{t}, this would provide a geometric interpretation of the classical Lie algebras analogous to his own achievement for the algebraic groups. In the PhD thesis of the second author \cite{r} this question was addressed in the case of simply-laced types and, in particular, a complete solution was obtained for the diagram $A_n$. In this paper we provide a uniform treatment of all simply-laced diagrams. Namely, we show that if $\E(L)$ is of simply-laced type, then $L$ is indeed a classical Lie algebra of the same type.

\begin{thm}\label{mainthm} Suppose that $L$ is a finite-dimensional simple Lie algebra generated by its extremal elements and contains no sandwiches. Furthermore, suppose that the extremal geometry $\E(L)$ is the root shadow space of a  building of type $A_n$ ($n\ge 2$), $D_n$ ($n \ge 4$), $E_6$, $E_7$ or $E_8$. Then $L$ is a classical Lie algebra of the same type as $\E(L)$.
\end{thm}

The outline of our proof is as follows. We fix an apartment $\Sigma$ in the building $\Delta$ related to $\E(L)$ and, using the known results about root shadow spaces, we deduce that $\Sigma$ gives a set of elements in $L$ that satisfy the familiar properties of a Chevalley basis. In particular, $L$ contains a classical subalgebra $L'$. From the Moufang condition of $\E(L)$, we are able to show that certain subgroups of $\aut(\E)$ associated with extremal elements have the same action on $\E(L)$ as the root subgroups of $\Delta$. In particular, using the standard theory of buildings, the group generated by such subgroups of $\aut(L)$ acts transitively on the apartments of $\E(L)$. We then deduce that the classical subalgebra $L'$ is in fact all of $L$.

The paper is organised in the following way. In Section \ref{sec:rss}, we describe the relationship between the root shadow space of an apartment, wherein points are described as arctic regions, and the corresponding root system. We survey the known results concerning root shadow spaces and construct a dictionary between root shadow spaces and two different interpretations of root systems. In Section \ref{sec:rfs}, we introduce root filtration spaces that are defined as particular filtrations on sets partitioned by five symmetric relations. We focus on an important example that emerges from the extremal elements of a Lie algebra $L$. In particular, the points of this  root filtration space are the $1$-spaces generated by extremal elements and we call it the extremal geometry $\E(L)$ of $L$. We then extend the dictionary constructed in Section \ref{sec:rfs} by adding the new language of extremal elements. Section \ref{sec:chevalg} is devoted to a 
characterisation of classical Lie algebras of simply-laced type
as algebras generated
by a so-called Chevalley spanning sets (a generalisation of Chevalley bases).
This characterisation is used in Section \ref{sec:subalg} to construct a classical subalgebra $L'$ of $L$ whenever $\E(L)$ is simply-laced. In Section \ref{sec:rootsubgroup}, we discuss the notion of a root subgroup of $\E(L)$ and show that for  each extremal element $x$, the map $\exp(x,t)$ produces such a group. Using this we deduce in the final section
that the classical subalgebra $L'$ coincides with $L$, finishing the proof
of the main result, Theorem \ref{mainthm}.

%--------------------------------------- Root Shadow Spaces --------------------------------------------------%

\section{Root shadow spaces and polar regions} \label{sec:rss}
In this section we explore the point-line spaces known as root shadow spaces and their relationship with polar regions of long roots. We view buildings as chamber systems and adopt the language of \cite{w2}. While this paper is mainly concerned with the simply-laced diagrams, the results we survey in this section also hold for non-simply-laced diagrams. Let $\Delta$ be an irreducible spherical building of type $X_n$ with corresponding root system $\Phi$. (Note that in the case where the building is of type $BC_n$, we only consider the type $B_n$). Let $I = \{1,\ldots,n\}$ be the index set of $\Delta$.

\begin{defn} Let $J$ be a fixed subset of $I$. A \emph{$J$-shadow} in $\Delta$ is an $(I\setminus J)$-residue. For each $j\in J$, a \emph{$j$-line} is the set of all $J$-shadows that contain chambers from a given $j$-panel. Let $\Pts$ be the set of $J$-shadows and $\Lines$ be the set of $j$-lines for $j\in J$. We say that $(\Pts, \Lines)$ is the \emph{$J$-shadow space} (or \emph{shadow space of type $X_{n,J}$}). If $J$ is the set $\{j\}$ then we write $X_{n,j}$ instead of $X_{n,J}$.
\end{defn}

In this paper we focus on a particular set $J$ defined as follows. Fix a fundamental system of simple roots $B = \{ \alpha_i \mid i\in I\}$ for $\Phi$ and denote the longest root with respect to $B$ by $\alpha$. Define the set $J\subset I$ consisting of all $i\in I$ such that $(\alpha,\alpha_i)\ne 0$. This corresponds to the set of nodes in the Dynkin diagram of $\Phi$ that are joined to the new node in its corresponding extended Dynkin diagram. We call $J$ the \emph{root set}. The root set is $J = \{1,n\}$ for $A_n$, $J=\{2\}$ for $(BC)_n, D_n, E_6$ and $G_2$, $J = \{ 1\}$ for $E_7$ and $F_4$, and $J = \{8\}$ for $E_8$. (We use the Bourbaki labelling of the Dynkin diagrams as in \cite{ci1,ci2,mw}). The root shadow space of type $X_{n,J}$, where $J$ is the root set, is simply called the \emph{root shadow space} of $\Delta$ and we denote it by $\rsh(\Delta)$. From now on $J$ is always a root set. Another description of the root shadow space of an irreducible building is in terms of its long root geometry (\cf \cite{Co1}).

Let $\Sigma$ be an apartment of $\Delta$. There is a natural one-to-one correspondence between  the half-apartments of $\Sigma$ and the roots of $\Phi$. The details of this can be found in Chapter 2 of \cite{w1}. It is common to simply identify the above two concepts with one another and, in particular, half-apartments are commonly referred to as roots. We will use the ``half-apartment'' terminology because we use both contexts in parallel in this section and Section \ref{sec:rootsubgroup}. This distinction in terminology will allow us to use the same notation for roots and half-apartments: say, if $\alpha$ is a root in $\Phi$, then we also use the same $\alpha$ to denote the corresponding half-apartment in $\Sigma$, and vice versa.

The arctic regions of half-apartments were first introduced in Chapter 6 of \cite{w2} and are later given a rigorous treatment in the paper \cite{mw}. The arctic region of a half-apartment $\alpha$ of $\Sigma$ is the set of chambers in $\alpha$ farthest away from the boundary $\partial\alpha$ of $\alpha$. More formally, we have the following definition.

\begin{defn} Let $\alpha$ be a half-apartment and $\Sigma$ be an apartment containing $\alpha$. The \emph{arctic region} of $\alpha$ is the set of chambers $x$ in $\alpha$ with the property that $\dist(x,\partial\alpha) \ge \dist(y,\partial\alpha)$ whenever $y$ is in $\alpha$ and $x$ and $y$ are incident. We denote the arctic region of $\alpha$ by $R_\alpha$.
\end{defn}

The arctic region $R_\alpha$ is a residue in $\Sigma$ (Lemma 6.9 of \cite{w2}). If $\alpha$ is a long root in $\Phi$, then $R_\alpha$ is a $J$-shadow of $\Sigma$ and every $J$-shadow arises in this way from a half-apartment corresponding to a long root of $\Phi$. Note that
$R_\alpha$ is defined in terms of $\alpha$ and $\partial\alpha$ and is therefore independent of the choice of apartment $\Sigma$. 
We now return to the building $\Delta$.

\begin{defn} The \emph{polar region} of a half-apartment $\alpha$ is the unique $J$-shadow of $\Delta$ containing the arctic region $R_\alpha$.
\end{defn}

Let $R$ be the polar region of a half-apartment $\alpha$ in $\Sigma$. Note that $\Sigma$ itself is a thin building and so we can consider the situation where $\Delta = \Sigma$. (Then polar regions and arctic regions coincide.) Hence we can consider the root shadow space $\rsh(\Sigma)$. Mapping every arctic region in $\Sigma$ to the (extended) polar region of $\Delta$ that contains it gives rise to an embedding of $\rsh(\Sigma)$ into $\rsh(\Delta)$. In fact, it will be convenient to identify $\rsh(\Sigma)$ with its image in $\rsh(\Delta)$.
We next recall a well-known dictionary between roots in $\Phi$, half-apartments in $\Sigma$ and $J$-shadows in $\rsh(\Sigma)$.

Fix an apartment $\Sigma$ and the corresponding root system $\Phi$. Let $\alpha$ and $\beta$ be long roots from $\Phi$ and $R_\alpha$ and $R_\beta$ be the corresponding arctic regions of $\Sigma$. Let $\theta_{\alpha\beta}$ denote the angle between $\alpha$ and $\beta$ in $\Phi$, let $d_{\alpha\beta}$ denote the distance between $R_\alpha$ and $R_\beta$ in the collinearity graph of $\rsh(\Sigma)$ and let $N_{\alpha\beta}$ denote the number of common neighbours of $R_\alpha$ and $R_\beta$ when they are non-collinear. Table \ref{tab:dict} describes the relationships between these parameters. Note that $-\alpha$ is the half-apartment complementary to $\alpha$ and that $\alpha+\beta$ is the half-apartment corresponding to the root $\alpha+\beta$ of $\Phi$.

\begin{table} 
\begin{center}
\begin{tabular}{|c|c|c|c|}
\hline
$\theta_{\alpha\beta}$ 		  & 	$d_{\alpha,\beta}$ 	& $N_{\alpha,\beta}$ 	& Half-apartments in $\Sigma$ \\
\hline
 $0$					 & 0 						&  						& $\alpha = \beta$\\
\hline
 $\pi/3$				 & 1						&  						& $R_\alpha\subset\beta$ \\
																		  &&& $R_\beta\subset\alpha$ \\
\hline
 $\pi/2$				 & 2						& $> 1$ 						& $R_\alpha\cap\beta \ne \emptyset$ \\
																		  &&& $R_\beta\cap\alpha \ne \emptyset$ \\
\hline
 $2\pi/3$			& 2							& 1 						& $R_\alpha\subset -\beta$ \\
												   && $(R_{\alpha+\beta})$      & $R_\beta\subset -\alpha$ \\
\hline
$\pi$				 & 3						& 0 						& $\alpha = -\beta$ \\
\hline

\end{tabular} 
\caption{Translation between the different interpretations of the long root geometry.} \label{tab:dict}
\end{center}
\end{table} 

The values of parameters $d_{\alpha,\beta}$ and $N_{\alpha,\beta}$ do not change if we substitute the arctic regions of $\Sigma$ with the corresponding polar regions of $\Delta$ and substitute $\rsh(\Sigma)$ with $\rsh(\Delta)$.   Note that for any two polar regions of $\Delta$ there is always an apartment $\Sigma$ that intersects both of them and so any two polar regions are in one of the five relations given in Table \ref{tab:dict}. All the details can be found in \cite{ci2,r}.

%--------------------------------------- End: Root Shadow Spaces --------------------------------------------------%

%--------------------------------------- Root Filtration Spaces --------------------------------------------------%

\section{Root filtration spaces and Lie algebras} \label{sec:rfs}

Root filtration spaces were introduced by Cohen and Ivanyos in \cite{ci1}. They were inspired by the examples of the so-called long root geometries which are in fact the same as the root shadow spaces introduced in Section \ref{sec:rss}.

A \emph{point-line space} (or simply a \emph{space}) is a pair $(\Pts, \Lines)$ where $\Pts$ is a set whose members are \emph{points} and $\Lines$ is a collection of subsets of $\Pts$, whose members are \emph{lines}. Two points of $\Pts$ are \emph{collinear} if they are contained in a common line in $\Lines$. For example, the root shadow spaces are examples of point-line spaces. A space is called a \emph{partial linear space} if every pair of distinct points lie on at most one common line, a \emph{singular space} if every two points lie on a common line, and a \emph{linear space} if it is both a singular and partial linear space. A subset $\X$ of $\Pts$ is a \emph{subspace} if it contains all the points of a line $l$ of $\Lines$ whenever $\X\cap l$ contains at least two points. Note that in a partial linear space every line is a subspace.  A line is \emph{thick} if it contains at least three points. A \emph{hyperplane} is a subspace that meets every line in at least one point. In particular, it either contains a line or intersects it in a unique point. The \emph{rank} of a linear space $\X$ is the length of a maximal chain of proper nontrivial subspaces. The rank of a trivial subspace is $0$ and the rank of a line is $1$. The \emph{singular rank} of a space is the supremum of the ranks of all maximal singular subspaces. 

The \emph{collinearity graph} of $(\Pts, \Lines)$ is
a graph whose vertex set is $\Pts$ and two distinct points are adjacent whenever they are collinear.

\begin{defn} \label{def:rfs} Let $(\E, \Flines)$ be a partial linear space and $\{\E_i\}_{-2\le i\le 2}$ be a set of symmetric relations on $\E$ partitioning $\E\times\E$. Let $\E_{\le j} = \cup_{i=-2}^j \E_i$. The space $(\E, \Flines)$ is a \emph{root filtration space with filtration $\{\E_i\}_{-2\le i\le 2}$} if the following hold.

\begin{itemize}
\item[{\rm{(A)}}] $\E_{-2}$ is equality between points in $\E$.
\item[{\rm{(B)}}] $\E_{-1}$ is collinearity between distinct points in $\E$.
\item[{\rm{(C)}}] There exists a map $\E_1 \map \E$, denoted by $(x,y)\mapto [x,y]$, such that $[x,y] \in \E_{\le i+j}(z)$ whenever $z\in\E_i(x)\cap\E_j(y)$.
\item[{\rm{(D)}}] If $(x,y)\in\E_2$, then $\E_{\le 0}(x)\cap\E_{\le -1}(y)$ is empty.
\item[{\rm{(E)}}] For each $x\in \E$, $\E_{\le -1}(x)$ and $\E_{\le 0}(x)$ are subspaces of $(\E, \Flines)$.
\item[{\rm{(F)}}] For each $x\in \E$, $\E_{\le 1}(x)$ is a hyperplane of $(\E, \Flines)$.
\end{itemize}

The space $(\E,\Flines)$ is a \emph{nondegenerate root filtration space} if the additional properties hold.

\begin{itemize}
\item[{\rm{(G)}}] For each $x\in\E$, $\E_2(x)$ is nonempty.
\item[{\rm{(H)}}] The collinearity graph $(\E, \E_{-1})$ is connected.
\end{itemize}
\end{defn}

Properties of root filtration spaces can be found in \cite{ci1} and \cite{ci2}. One important consequence of such properties is the following result.

\begin{lem} \label{points-lines} Let $(\E, \Flines)$ be a nondegenerate root filtration space. Then its defining relations can be characterised by the collinearity graph $(\E, \E_{-1})$ in the following way:

\begin{itemize}
\item[{\rm{(-2)}}] $(x,y)\in\E_{-2}$ if and only if $x=y$;
\item[{\rm{(-1)}}] $(x,y)\in\E_{-1}$ if and only if $x$ and $y$ are distinct collinear points;
\item[{\rm{(0)}}] $(x,y)\in\E_0$ if and only if $x$ and $y$ are non-collinear and have at least two common neighbours;
\item[{\rm{(1)}}] $(x,y)\in\E_1$ if and only if $x$ and $y$ are non-collinear and have a unique common neighbour;
\item[{\rm{(2)}}] $(x,y)\in\E_2$ if and only if $x$ and $y$ are non-collinear and have no common neighbours.
\end{itemize}

Furthermore, pairs of points in $\E_{-2}$, $\E_{-1}$, $\E_0 \cup \E_1$, and $\E_2$ have distances between them equal to $0$, $1$, $2$ and $3$, respectively. \qed
\end{lem}

Essentially, the structure of a nondegenerate root filtration space can be recovered from the relation $\E_{-1}$.

We often drop the $\Flines$ and simply refer to the root filtration space as $\E$. Any nondegenerate root filtration space $\E$ with finite singular rank can be identified with a root shadow space of a building. This important result was proven, in a slightly different language, by Kasikova and Shult in \cite{ks} when $\E$ has polar rank at least $3$ and by Cohen and Ivanyos in \cite{ci2} when $\E$ has polar rank $2$. The proof entails the reconstruction of the building from $\E$.

\begin{thm} \label{mainresult} A nondegenerate root filtration space with finite singular rank is isomorphic to the root shadow space of type $A_{n,\{1,n\}}$, $BC_{n,2}$, $D_{n,2}$, $E_{6,2}$, $E_{7,1}$, $E_{8,8}$, $F_{4,1}$, or $G_{2,2}$. \qed
\end{thm}

Note that, formally speaking, Cohen and Ivanyos in \cite{ci2} only claim isomorphism of the point-line spaces. Lemma \ref{points-lines} shows that the structure of the root filtration space can be uniquely recovered from its point-line space. Hence we do indeed get the isomorphism that preserves the five relations $\E_i$ of the root filtration space as a consequence of \cite{ci2}. (Undoubtedly, Cohen and Ivanyos were aware of this.) 

We now describe the root filtration space associated to the so-called extremal elements of a Lie algebra. Firstly, we need the precise definition of an extremal element which works for any field characteristic.

\begin{defn} \label{extremal elements} Let $L$ be a Lie algebra over a field $F$. An element $x\in L$ is said to be \emph{extremal} if there is a map $g_x : L\map F$ such that, for all $y$ and $z$ in $L$, we have

\begin{equation} \label{extremal} [x,[x,y]] = 2g_x(y)x,
\end{equation}
\noindent and
\begin{equation} \label{exchar2a} [[x,y],[x,z]] = g_x([y,z])x + g_x(z)[x,y] - g_x(y)[x,z]
\end{equation}
\end{defn}

If $F$ does not have characteristic $2$, then (\ref{extremal}) implies (\ref{exchar2a}). Hence, we end up with the familiar simple definition as in the introduction. Furthermore, in this case $g_x$ is clearly unique. If $F$ has characteristic $2$, $g_x$ is not necessarily unique. For example, if $L$ is abelian, then $g_x(0) = 0$ is the only restriction. 

An extremal element $x\in L$ is a \emph{sandwich} whenever $g_x$ can be taken to be zero. By convention, we always take $g_x$ to be zero for sandwiches. 

Under this convention, we can state the following result, which is essentially Proposition 20 from \cite{ci1}.

\begin{lem} \label{ext=kill}
Suppose that $L$ is generated by its extremal elements. Then there is a unique bilinear form $g:L\times L\map F$ such that $g(x,y) = g_x(y)$ for every extremal element $x$. The form $g$ is symmetric and invariant under the Lie product. Namely, $g_x(y) = g_y(x)$ and $g_x([y,z]) = g_{[x,y]}(z)$ for all $x,y,z\in L$. Furthermore, $g_x(y) = 0$ whenever $[x,y]=0$.
\end{lem} 

In the remainder of the paper, the Lie algebras we consider do not contain sandwich elements. Thus, extremal elements are considered to be exclusively non-sandwich extremal elements. In particular, the zero element is not considered extremal. Let $E$ be the set of extremal elements of $L$. Define five symmetric binary relations $E_i$ for $i\in\{-2,-1,0,1,2\}$ on $E$ as follows:

\begin{itemize} \label{therelns}
\item[{\rm{(-2)}}] $(x,y)\in E_{-2}$ if and only if $x$ and $y$ are linearly dependent;
\item[{\rm{(-1})}] $(x,y)\in E_{-1}$ if and only if $x$ and $y$ are linearly independent,
$[x,y]=0$, and $\lambda x + \mu y \in E$ for all $(\lambda, \mu)\in F^2$, $(\lambda,\mu)\ne (0,0)$;
\item[{\rm{(0)}}] $(x,y)\in E_0$ if and only if $[x,y]=0$ and $(x,y)$ is not in $E_{-2}\cup E_{-1}$;
\item[{\rm{(1)}}] $(x,y)\in E_1$ if and only if $[x,y]\ne 0$, but $g_x(y) = 0$;
\item[{\rm{(2)}}] $(x,y)\in E_2$ if and only if $g_x(y) \ne 0$.
\end{itemize}

\medskip
Note that $g_x(y) = 0$ whenever $(x,y)\in \cup_{j=-2}^1 E_j$ and also $[x,y]\ne 0$ for all $(x,y)\in\E_1\cup\E_2$. 

Let $\E$ be the set of projective points spanned by the extremal elements of $L$, that is, $\E = \{ Fx \mid x\in E\}$. Let $\Flines$ be the set of projective lines $Fx+Fy$ for $(x,y)\in E_{-1}$. We call $(\E$, $\Flines)$ the \emph{extremal geometry} of the Lie algebra $L$ and denote it by $\E(L)$ to indicate it is derived from $L$. We call the points in $\E$ \emph{extremal points} (or simply \emph{points}). By the definition of $E_{-1}$, the unique line in $\Flines$ containing the two collinear points $Fx$ and $Fy$ is $Fx+Fy$. Hence $\E(L)$ is a partial linear space. The symmetric relations $\{\E_i\}_{i=-2}^2$ correspond to $\{E_i\}_{i=-2}^2$ in the natural way, namely, $(Fx, Fy)\in\E_i$ if and only if $(x,y)\in E_i$. In particular, $\{\E_i\}_{i=-2}^2$ are five disjoint symmetric relations on $\E$ that partition $\E\times\E$, where $\E_{-2}$ is equality and $\E_{-1}$ is collinearity in $\E(L)$. The next result is a corollary of Theorem 28 of \cite{ci1}.

\begin{thm} \label{nondegext} Let $L$ be a finite-dimensional simple Lie algebra over $F$ generated by its extremal elements. Then $\E(L)$ is either a nondegenerate root filtration space with finite singular rank whose lines are all thick or a root filtration space without lines. \qed
\end{thm}

Note that in the case where $\E(L)$ has lines it follows from Theorem \ref{mainresult} that $\E(L)$ is isomorphic to $\rsh(\Delta)$ for an irreducible spherical building $\Delta$.

%--------------------------------------- End: Root Filtration Spaces --------------------------------------------------%

%--------------------------------------- Classical algebras ---------------------------------------%

\section{Classical algebras} \label{sec:chevalg}

In this section we provide a recognition result for classical algebras of simply-laced type. 

Let $\Phi$ be an irreducible root system with a simply-laced diagram $X_n$ and let $(|)$ be the standard inner product on the Euclidean space spanned by $\Phi$.

\begin{defn} \label{spanning set}
We say that a Lie algebra $L$ over a field $F$ is \emph{of type $X_n$} if it possesses a spanning set consisting of elements $x_\alpha$ and $h_\alpha$, $\alpha\in\Phi$, which satisfy the following, for $\alpha,\beta\in\Phi$:
\begin{itemize}
\item[{\rm{(i)}}] $[h_\alpha, h_\beta] = 0$;
\item[{\rm{(ii)}}] $[h_\alpha, x_\beta] = (\alpha|\beta) x_\beta$;
\item[{\rm{(iii)}}] if $\alpha+\beta=0$ then $[x_\alpha,x_\beta]=h_\alpha$; and
\item[{\rm{(iv)}}] if $\alpha+\beta\neq 0$ then $[x_\alpha, x_\beta] = A_{\alpha,\beta} x_{\alpha+\beta}$ where  $A_{\alpha,\beta}\neq 0$ if and only if $\alpha+\beta$ is a root.
\end{itemize}
We call such a set a \emph{Chevalley spanning set}. 
\end{defn}

From (iii), $h_\alpha=[x_\alpha,x_{-\alpha}]$, and this immediately gives the following. 

\begin{lem}
For $\alpha\in\Phi$, we have $h_{-\alpha}=-h_\alpha$. \qed
\end{lem}

Clearly, we can scale $x_\alpha$ and $x_{-\alpha}$ by $c$ and $c^{-1}$, respectively, for any $0\ne c\in F$, and this operation again gives a Chevalley spanning set in $L$. Hence we have a lot of freedom in choosing the structure constants $A_{\alpha,\beta}$. We will see later that we can ask for all nonzero $A_{\alpha,\beta}$ to be equal to $\pm 1$ and, furthermore, the signs can also be controlled. However, first we record some general properties. 

Note that $A_{\alpha,\beta}\neq 0$ if and only if $\alpha+\beta\in\Phi$ and so if and only if the angle between $\alpha$ and $\beta$ is $2\pi/3$.

\begin{lem} \label{circle}
Suppose that the roots $\alpha$ and $\beta$ form an angle of $2\pi/3$. Then $A_{\alpha,\beta}=-A_{\beta,\alpha}$ and $A_{\alpha,\beta}=-A_{\alpha+\beta,-\alpha}^{-1}=
A_{\beta,-\alpha-\beta}=-A_{-\alpha,-\beta}^{-1}=A_{-\alpha-\beta,\alpha}=-A_{-\beta,\alpha+\beta}^{-1}$.
\end{lem}

\pf Only the second claim requires proving. Note that $(\alpha|\beta)=-1$. Hence, by (ii),  $-x_\beta=[h_\alpha,x_\beta]=[[x_\alpha,x_{-\alpha}],x_\beta]=[[x_\alpha,x_\beta],x_{-\alpha}]$. (We used that $-\alpha+\beta$ is not a root and hence $[x_{-\alpha},x_\beta]=0$.) Finally, since $[x_\alpha,x_\beta]=A_{\alpha,\beta}x_{\alpha+\beta}$, we obtain $-x_\beta=[A_{\alpha,\beta}x_{\alpha+\beta},x_{-\alpha}]$, and so $[x_{\alpha+\beta},x_{-\alpha}]=-A_{\alpha,\beta}^{-1}x_\beta$, proving $A_{\alpha+\beta,-\alpha}=-A_{\alpha,\beta}^{-1}$, or equivalently, $A_{\alpha,\beta}=-A_{\alpha+\beta,-\alpha}^{-1}$. Applying this with $\alpha+\beta$ and $-\alpha$ in place of $\alpha$ and $\beta$, we get the second equality. Repeating this process around the $A_2$ root system in the plane spanned by $\alpha$ and $\beta$ gives the remaining equalities. \qed

\begin{lem} \label{sum}
For $\alpha,\beta\in\Phi$, if $\alpha+\beta\in\Phi$, then $h_{\alpha+\beta}=h_\alpha+h_\beta$. 
\end{lem}

\pf From (iv) we have $A_{\alpha,\beta}x_{\alpha+\beta}=[x_\alpha,x_\beta]$. Hence,
\begin{align*} 
h_{\alpha+\beta}&=[x_{\alpha+\beta},x_{-\alpha-\beta}]\\
 &= A_{\alpha,\beta}^{-1}[[x_\alpha,x_\beta],x_{-\alpha-\beta}]\\
 &=A_{\alpha,\beta}^{-1}([[x_\alpha,x_{-\alpha-\beta}],x_\beta]+[x_\alpha,[x_\beta,x_{-\alpha-\beta}]])\\
 &=A_{\alpha,\beta}^{-1}(A_{\alpha,-\alpha-\beta}[x_{-\beta},x_\beta]+      A_{\beta,-\alpha-\beta}[x_\alpha,x_{-\alpha}]).
\end{align*}
By Lemma \ref{circle}, $A_{\alpha,\beta}^{-1}A_{\alpha,-\alpha-\beta}=-1$ and $A_{\alpha,\beta}^{-1}A_{\beta,-\alpha-\beta}=1$.  Therefore, $h_{\alpha+\beta}=-[x_{-\beta},x_\beta]+[x_\alpha,x_{-\alpha}]= h_\beta+h_\alpha$, and the claim follows. \qed

\medskip
Let $I = \{ 1, \ldots, n \}$ be an index set and $\{\alpha_i \mid i\in I\}$ be a fundamental system of simple roots of $\Phi$. Set $h_i := h_{\alpha_i}$ for $i\in I$. Taking into account that every root is an integral linear combination of simple roots, we have the following result. 

\begin{cor} \label{dependence}
For $\alpha\in\Phi$, we have $h_\alpha=\sum_{i\in I}c_ih_i$, where the integral coefficients $c_i$ are defined by $\alpha=\sum_{i\in I}c_i\alpha_i$.
\end{cor}

\pf Clearly, we can assume that $\alpha$ is a positive root. Since every such root is a sum of two positive roots of lower height, the claim follows by induction. \qed

\begin{cor} \label{smaller}
The set of all $x_\alpha$, $\alpha\in\Phi$, and $h_i$, $i\in I$, spans $L$. \qed
\end{cor}

We now turn to scaling the vectors $x_\alpha$ in order to obtain the simplest possible structure constants $A_{\alpha,\beta}$. First of all, for each positive non-simple root $\alpha$, we select a particular pair $(\alpha_1,\alpha_2)$ of positive roots satisfying $\alpha=\alpha_1+\alpha_2$. We will refer to $(\alpha_1,\alpha_2)$ as the \emph{defining pair} of $\alpha$. Clearly, there are many ways to choose defining pairs, we just make our selection in some way and fix it.

We will now scale our vectors $x_\alpha$ inductively as follows. (Recall that whenever we scale $x_\alpha$ for a positive root $\alpha$ by a nonzero scalar $c\in F$, by convention, we will also scale $x_{-\alpha}$ by $c^{-1}$ so that $h_\alpha=[x_\alpha,x_{-\alpha}]$ remains unchanged.) If $\alpha=\alpha_i$ is a simple root (height 1), then we leave $x_\alpha$ unchanged. Any other positive root $\alpha$ has a defining pair $(\alpha_1,\alpha_2)$, and we scale $x_\alpha$ so that $x_\alpha=[x_{\alpha_1},x_{\alpha_2}]$. Note that $\alpha_1$ and $\alpha_2$ have lower heights than $\alpha$, so the inductive scaling procedure is well-defined. The consequence of this operation is that we now have that $A_{\alpha_1,\alpha_2}=1$ for all defining pairs $(\alpha_1,\alpha_2)$.

\begin{lem} \label{plus minus one}
Suppose that $A_{\gamma_1,\gamma_2}=1$ whenever $(\gamma_1,\gamma_2)$ is a defining pair. Then all nonzero structure constants $A_{\alpha,\beta}$ are equal to $\pm 1$ and the signs are fully determined by the choice of defining pairs in $\Phi$.
\end{lem}

\pf The $A_2$ system in the plane spanned by $\alpha$ and $\beta$ contains a single pair of positive roots $\alpha'$ and $\beta'$ at the angle $2\pi/3$. Furthermore, Lemma \ref{circle} shows that $A_{\alpha,\beta}=\pm A_{\alpha',\beta'}^{\pm 1}$. So the claim holds for $\alpha$ and $\beta$ if and only if it holds for $\alpha'$ and $\beta'$. Hence we just need to consider the case where both $\alpha$ and $\beta$ are positive.

Let $\gamma=\alpha+\beta$. We do induction on the height of $\gamma$. Let $(\gamma_1,\gamma_2)$ be the defining pair for $\gamma$. By assumption, $x_\gamma=[x_{\gamma_1},x_{\gamma_2}]$, since $A_{\gamma_1,\gamma_2}=1$. As $\alpha+\beta=\gamma_1+\gamma_2$, the root subsystem $\Psi$ spanned by $\{ \alpha, \beta, \gamma_1, \gamma_2 \}$ has rank at most $3$ and so is either of type $A_2$ or $A_3$. If $\Psi$ has type $A_2$, then $\{\alpha,\beta\} = \{\gamma_1, \gamma_2\}$ and thus our claim holds. Suppose now that $\Psi$ has type $A_3$. Note that $\Psi$ inherits the partition into positive and negative roots from $\Phi$. Since $\gamma$ has two different decompositions as a sum of positive roots in $\Psi$, it is clear that $\gamma$ is the highest positive root of $\Psi$. Also, in each of the two pairs $(\alpha,\beta)$ and $(\gamma_1, \gamma_2)$, one of the roots is simple and the other has height $2$. This gives four possibilities which are all similar and so we just consider one of them. Let us suppose that $\alpha$ and $\gamma_2$ are simple roots. Let $\kappa$ be the remaining third simple root in $\Psi$. Note that $\alpha + \kappa + \gamma_2 = \gamma = \alpha+\beta = \gamma_1 + \gamma_2$ and so $\beta = \kappa + \gamma_2$ and $\gamma_1 = \alpha + \gamma$. Thus $\alpha$ and $\gamma_2$ are perpendicular and both form an angle of $2\pi/3$ with $\kappa$.
Finally, we have
\[
	[x_\alpha, x_\beta] = [x_\alpha, A_{\kappa,\gamma_2}^{-1}[x_\kappa, x_{\gamma_2}]] = A_{\kappa,\gamma_2}^{-1}[x_\alpha, [x_\kappa,x_{\gamma_2}]].
\]
Since $\alpha$ and $\gamma_2$ are perpendicular, we have $[x_\alpha,x_{\gamma_2}]=0$ and $[x_\alpha, [x_\kappa,x_{\gamma_2}]]=[[x_\alpha,x_\kappa],x_{\gamma_2}]$. Therefore,
\[ [x_\alpha, x_\beta]=A_{\kappa,\gamma_2}^{-1}[[x_\alpha,x_\kappa],x_{\gamma_2}]=A_{\kappa,\gamma_2}^{-1} A_{\alpha,\kappa}[x_{\gamma_1},x_{\gamma_2}]=A_{\kappa,\gamma_2}^{-1}A_{\alpha,\kappa}x_\gamma.
\]
Thus $A_{\alpha,\beta}=A_{\kappa,\gamma_2}^{-1}A_{\alpha,\kappa}$. By induction, the two factors on the right are $\pm 1$ and are determined by the choice of defining pairs. Therefore, the same is true for $A_{\alpha, \beta}$. \qed

\medskip
Consider now the simple complex algebra $\g$ corresponding to the simply-laced root system $\Phi$. It is well known that $\g$ has a Chevalley spanning set. Indeed, it can be constructed from a Cartan decomposition of $\g$. Additionally, the smaller spanning set, $\{x_\alpha,h_i\mid\alpha\in\Phi,i\in I\}$ is a basis of $\g$. When all structure constants $A_{\alpha,\beta}$ are equal to $\pm 1$, a basis of this sort is known as a \emph{Chevalley basis} of $\g$. It follows from Corollary \ref{dependence} that all structure constants for such a basis are integer, which means that the integral span $\g_\Z$ of the Chevalley basis is closed under Lie product. For an arbitrary field $F$, $\hat\g=\g_F=\g_\Z\otimes_\Z F$ is a Lie algebra over $F$. We call $\hat \g$ the \emph{Chevalley algebra of type $X_n$ over $F$.} Note the Chevalley algebra is well-defined, namely, it does not depend on the choice of the Chevalley basis. Indeed, given two Chevalley bases, we can assume, after scaling $x_\alpha$ and $x_{-\alpha}$ by $\pm 1$ as necessary, that both bases have $A_{\alpha_1,\alpha_2}=1$ for all defining pairs $(\alpha_1,\alpha_2)$. It now follows from Corollary \ref{dependence} and Lemma \ref{plus minus one} that the two bases have identical structure constants, which means that they are conjugate by an automorphism of $\g$. The above scaling by $\pm 1$ does not affect the integral span of the basis, hence the algebras $\g_\Z$ corresponding to different Chevalley bases are conjugate. Manifestly, this means that the Chevalley algebras $\hat\g=\g_\Z\otimes_\Z F$ coming from different Chevalley bases of $\g$ must be isomorphic.

If we set $\hat x_\alpha:=x_\alpha\otimes 1_F$ and $\hat h_\alpha:=h_\alpha\otimes 1_F$, then the elements $\hat x_\alpha$ and $\hat h_\alpha$, $\alpha\in\Phi$, form a Chevalley spanning set for $\hat\g$, while the elements $\hat x_\alpha$, $\alpha\in\Phi$, together with $\hat h_i$, $i\in I$, form a basis of $\hat \g$. Recall that we defined classical algebras as nonzero factor algebras of Chevalley algebras. The image in the factor algebra of a Chevalley spanning set is again a Chevalley spanning set. All the required properties are immediate as long as no element from the spanning set becomes zero. If the image of $x_\alpha$ is zero, then also the image of $h_\alpha=[x_\alpha,x_{-\alpha}]$ is zero. This in turn implies that all $x_\beta$ map to zero, where the angle between $\alpha$ and $\beta$ equals $\pi/3$ or $2\pi/3$. This is because $[h_\alpha,x_\beta]=\pm x_\beta$. By connectivity we then have that all elements of the spanning set map to zero, which is a contradiction since the factor algebra is nonzero. We have shown that all classical algebras contain a Chevalley spanning set and hence they are all algebras of type $X_n$ as defined in the beginning of the section. 

We now come to the main result of this section.

\begin{prop} \label{recognition}
Let $\Phi$ be an irreducible root system of  simply-laced type $X_n$.
A Lie algebra $L$ has a Chevalley spanning set of type $X_n$, if and only if it is classical of type $X_n$.
\end{prop}

\pf We just need to show that every Lie algebra $L$ over a field $F$ with a Chevalley spanning set is isomorphic to a factor algebra of the Chevalley algebra $\hat\g=\g_F$ of the same type $X_n$. Choose Chevalley spanning sets $\{\hat x_\alpha,\hat h_\alpha\mid\alpha\in\Phi\}$ and $\{x_\alpha,h_\alpha\mid\alpha\in\Phi\}$ in $\hat\g$ and $L$ respectively.  Without loss of generality, both sets have the structure constants $A_{\alpha_1,\alpha_2}$ equal $1$ for all defining pairs $(\alpha_1,\alpha_2)$. Then by Lemma \ref{plus minus one} these two sets have identical values of $A_{\alpha,\beta}$ for all $\alpha$ and $\beta$ with $\alpha+\beta\ne 0$.

By the above, the set $\{\hat x_\alpha,\hat h_i\mid\alpha\in\Phi,i\in I\}$ is a basis of $\hat\g$, where $\hat h_i=\hat h_{\alpha_i}$. This means that there is a linear mapping $\phi:\hat\g\to L$ sending every $\hat x_\alpha$ to $x_\alpha$ and every $\hat h_i$ to $h_i=h_{\alpha_i}$. Corollary \ref{dependence} and the fact the constants $A_{\alpha,\beta}$ are the same for the two spanning sets imply that $\phi$ is an algebra homomorphism. Corollary \ref{smaller} implies that $\phi$ is onto, and so $L$ is isomorphic to a factor algebra of $\hat\g$. \qed

%--------------------------------------- End: Classical algebras ---------------------------%

%--------------------------------------- The Chevalley Algebra --------------------------------------------------%

\section{Constructing a classical subalgebra in $L$}
\label{sec:subalg}

Suppose $L$ is a finite-dimensional simple Lie algebra over the field $F$ that is generated by its extremal elements $E$. Furthermore suppose that its extremal geometry $\E = \E(L)$ is isomorphic to the root shadow space $\rsh(\Delta)$ of a building $\Delta$ of simply-laced type $X_n$ with index set $I = \{1, \ldots, n\}$, where $n\ge 2$. In this section we show that $L$ contains a classical subalgebra $L'$. 

Let $\Sigma$ be a fixed apartment of $\Delta$. Let $\Phi$ be the corresponding root system and let $\{\alpha_i \mid i\in I\}$ be a fundamental system of simple roots of $\Phi$. Recall that every root $\alpha$ in $\Phi$ corresponds to a half-apartment in $\Sigma$ (also denoted $\alpha$) which in turn corresponds to its arctic region $R_\alpha$. The arctic region $R_\alpha$ uniquely extends to a polar region and that is a $J$-shadow of $\Delta$, which corresponds to a point in the extremal geometry $\E$. In particular, there is a natural injection from $\Phi$ into $\E$. In what follows, we let $x_\alpha$ be an extremal element representing the point of $\E$ corresponding to $\alpha\in\Phi$. Let $\theta_{\alpha\beta}$ be the angle between roots $\alpha$ and $\beta$. 

\begin{lem} \label{anglesext} For roots $\alpha$ and $\beta$ in $\Phi$, we have the following;
\begin{itemize}
\item[{\rm{(i)}}] if $\theta_{\alpha\beta}\in\{0,\pi/3,\pi/2\}$, then $[x_\alpha,x_\beta]=0$;
\item[{\rm{(ii)}}] if $\theta_{\alpha\beta} = 2\pi/3$, then $F[x_\alpha, x_\beta]= Fx_{\alpha+\beta} \ne 0$;
\item[{\rm{(iii)}}] $g_{x_\alpha}(x_\beta)\ne 0$ if and only if $\theta_{\alpha\beta}=\pi$.
\end{itemize}
\end{lem}

\pf  The result follows from the definition of the extremal geometry $\E$, Theorem \ref{nondegext}, and Table \ref{tab:dict}. \qed

\medskip
Let $L'$ be the subalgebra of $L$ generated by $\{x_\alpha \mid \alpha\in\Phi\}$. The goal is to construct a Chevalley spanning set in $L'$. So far we have selected vectors $x_\alpha$ arbitrarily within the corresponding $1$-spaces of $L$. Let us now impose some restrictions. Namely, in view of Lemma \ref{anglesext} (iii), for each positive $\alpha\in\Phi$, we can scale $x_{-\alpha}$ so that $g_{x_\alpha}(x_{-\alpha})=-1$. Recall that $g_x(y)=g_y(x)$ and so $g_{x_{-\alpha}}(x_\alpha)=-1$. Finally, for each $\alpha$, positive or negative, we set $h_\alpha = [x_{\alpha}, x_{-\alpha}]$ and we note that $h_{-\alpha}=-h_\alpha$. 

We start with the products $[h_\alpha,x_\beta]$.

\begin{lem} \label{products with h}
For arbitrary roots $\alpha$ and $\beta$, we have $[h_\alpha,x_\beta]=(\alpha|\beta)x_\beta$.
\end{lem}

\pf If the angle between $\alpha$ and $\beta$ is zero, then $\alpha=\beta$ and so $[h_\alpha,x_\beta]= [[x_\alpha,x_{-\alpha}],x_\alpha]=[x_\alpha,[x_{-\alpha},x_\alpha]]=-[x_\alpha,[x_\alpha,x_{-\alpha}]]= -2g_{x_\alpha}(x_{-\alpha})x_\alpha=2x_\alpha=(\alpha|\alpha)x_\alpha$ by the choice of scaling above.

Next consider the angle $\pi/3$. Take $x=x_\alpha$, $y=x_{-\alpha}$, and $z=x_{\beta-\alpha}$ in the equation (\ref{exchar2a}) of Definition \ref{extremal elements}. This gives $[[x_\alpha,x_{-\alpha}],[x_\alpha,x_{\beta-\alpha}]]= g_{x_\alpha}([x_{-\alpha},x_{\beta-\alpha}])x_\alpha+g_{x_\alpha}(x_{\beta-\alpha})[x_\alpha,x_{-\alpha}]- g_{x_\alpha}(x_{-\alpha})[x_\alpha,x_{\beta-\alpha}]$. By noting that $[x_{-\alpha},x_{\beta-\alpha}]=0$ (because the angle here is $\pi/3$) and $g_{x_\alpha}(x_{\beta-\alpha})=0$ (because the angle between $\alpha$ and $\beta-\alpha$ is not $\pi$), we see that the above simplifies to $[h_\alpha,[x_\alpha,x_{\beta-\alpha}]]=  -g_{x_\alpha}(x_{-\alpha})[x_\alpha,x_{\beta-\alpha}]=[x_\alpha,x_{\beta-\alpha}]$. Finally, $[x_\alpha,x_{\beta-\alpha}]=ax_\beta$ for a nonzero scalar $a\in F$ in view of Lemma \ref{anglesext} (ii). Therefore we have $[h_\alpha,ax_\beta]=ax_\beta$, which clearly implies that $[h_\alpha,x_\beta]= x_\beta=(\alpha|\beta)x_\beta$.

If the angle between $\alpha$ and $\beta$ is $\pi/2$, then $[h_\alpha,x_\beta]=[[x_\alpha,x_{-\alpha}],x_\beta]= [[x_\alpha,x_\beta],x_{-\alpha}]+[x_\alpha,[x_{-\alpha},x_\beta]]=[0,x_{-\alpha}]+[x_\alpha,0]=0= (\alpha|\beta)x_\beta$.

Finally, if the angle is greater than $\pi/2$, then $[h_\alpha,x_\beta]=[-h_{-\alpha},x_\beta]= -(-\alpha|\beta)x_\beta=(\alpha|\beta)x_\beta$ since the angle between $-\alpha$ and $\beta$ is less than $\pi/2$.
\qed

\medskip
Next we consider the products $[h_\alpha,h_\beta]$.

\begin{lem} \label{Hcomm} 
For $\alpha,\beta\in\Phi$, we have $[h_{\alpha},h_{\beta}]=0$. 
\end{lem}

\pf Using the Jacobi identity, we write $[h_\alpha,h_\beta]=[h_\alpha,[x_\beta, x_{-\beta}]]=[[h_\alpha, x_\beta], x_{-\beta}] + [x_\beta,[h_\alpha, x_{-\beta}]]$. Applying Lemma \ref{products with h} yields $[(\alpha|\beta)x_\beta,x_{-\beta}]+[x_\beta,(\alpha|-\beta)x_{-\beta}]=((\alpha|\beta)+(\alpha|-\beta))h_\beta= (\alpha|\beta-\beta)h_\beta=0$.
\qed

\medskip
It follows from Lemmas \ref{anglesext}, \ref{products with h}, and \ref{Hcomm}, that the subspace spanned by $ \{ x_\alpha, h_\alpha \mid \alpha\in\Phi\}$ is closed under Lie product, and so it coincides with $L'$. This gives us a spanning set in $L'$ and, manifestly, this set satisfies all properties of a Chevalley spanning set. Hence Proposition \ref{recognition} yields the following result.

\begin{prop} \label{LisChev} 
Suppose $L$ is a simple Lie algebra over $F$ generated by extremal elements and suppose the extremal geometry $\E$ of $L$ comes from a building $\Delta$ of a simply-laced type $X_n$. Then the elements $x_\alpha$ derived from an arbitrary apartment $\Sigma$ of $\Delta$ generate a classical subalgebra $L'$ in $L$ of the same type $X_n$. \qed
\end{prop}

%--------------------------------------- End: The Chevalley Algebra --------------------------------------------------%

%--------------------------------------- The Roots Subgroups --------------------------------------------------%

\section{Root subgroups}
\label{sec:rootsubgroup}
 
In the final two sections we will study automorphisms of geometrical structures. If $X$ is a geometry (or a building), 
$g\in\aut(X)$ and $x\in X$, then we denote $x^g$ to be the image of $x$ under $g$.

Throughout this section $\E$ is a nondegenerate root filtration space. By an automorphism of $\E$ we mean a 
permutation of $\E$ that preserves the collinearity relation $\E_{-1}$. Recall that every other relation $\E_i$ 
is recoverable from $\E_{-1}$, which means that an automorphism preserves every $\E_i$. As usual, $\aut(\E)$ 
denotes the full group of automorphisms of $\E$. Additionally, we
will frequently switch between the language of the point-line space $\E$ and its collinearity graph (which determine one another) for convenience. The notions of \emph{path, neighbour} and \emph{adjacency} are reserved for the graph and \emph{lines} and \emph{collinearity} are reserved for the space.   

For a point $x\in\E$ and a line $l$, we say that $l$ is a \emph{cutting line} with respect to $x$ if $l$ contains 
points from $\E_{\leq -1}(x)$ and $\E_1(x)$. Note that since $\E_{\le 0}(x)$ is a subspace of $\E$, the cutting 
line $l$ contains exactly one point from $\E_{\le 0}(x)$ and clearly this point must lie in $\E_{-1}(x)$. Hence 
$l$ ``cuts from'' $\E_1(x)$ straight into $\E_{-1}(x)$ bypassing $\E_0(x)$. This explains the terminology.
Since $\E_{\le 1}(x)$ is also a subspace (in fact, a hyperplane) and since $l$ meets $\E_{\le 1}(x)$ in at least 
two points, we have that $l$ is fully contained in $\E_{\le 1}(x)$, that is, all the remaining points of $l$ 
are in $\E_1(x)$. Note that every point $y\in\E_1(x)$ lies in exactly one cutting line, namely the line through 
$y$ and $[x,y]$. This is because $[x,y]$ is the only common neighbour of $x$ and $y$. Hence the cutting lines 
with respect to $x$ lead to a partition of $\E_1(x)$.

\begin{defn} \label{point group}
For a point $x\in \E$, we define the \emph{point group} $Y(x)$ to be the subgroup of $\aut(\E)$ fixing $\E_{\le 0}(x)$ 
pointwise and stabilising every cutting line with respect to $x$ setwise.
\end{defn}

Clearly, $Y(x)$ acts on every cutting line $l$ and on $\E_1(x)$. 
We aim to show that the latter action is semiregular, that is, the 
stabiliser in $Y(x)$ of every $y\in\E_1(x)$ is trivial. We first need some technical lemmas.

We call two lines in $\E$ \emph{opposite} if for every point $x$ on one of them there is a point $y$ on the other 
line such that $(x,y)\in\E_2$. Suppose $l$ and $l'$ are opposite. Since $\E_{\leq 1}(z)$ is a hyperplane, it follows 
that for every $z$ on $l$ there exists a unique $z'$ on $l'$ such that $z'\in\E_{\leq 1}(z)$ and all other points of 
$l'$ lie in $\E_2(z)$. Note that by the property (D) of Definition \ref{def:rfs} of a root filtration space we have 
that $(z,z')\in\E_1$ and that the correspondence $z\to z'$ is a bijection between $l$ and $l'$.

\begin{lem} \label{opposite}
The following hold.
\begin{enumerate}
\item[\rm (i)] For distinct collinear points $x$ and $y$, there exists a point $z$ collinear with $y$ such that 
$(x,z)\in\E_1$.
\item[\rm (ii)] If $xuvz$ is a path in the collinearity graph of $\E$ such that $(x,v)$ and $(u,z)$ are in $\E_1$, then 
$(x,z)\in\E_2$.
\item[\rm(iii)] If $xuvwy$ is a path in the collinearity graph of $\E$ such that $(x,v)$, $(u,w)$ and $(v,y)$ are all 
contained in $\E_1$, then the lines $xu$ and $wy$ are opposite.
\end{enumerate}
\end{lem}

\pf The claim (i) follows from Lemma 4 of \cite{ci1} since $\E$ is nondegenerate. Claim (ii) is Lemma 1 (v) from 
\cite{ci1}. For the final claim, take $z$ on the line $xu$. If $z=u$ then choosing $z'=y$ gives $(z,z')\in\E_2$. On 
the other hand, if $z\ne u$ then taking $z'=w$ yields the same result. Indeed, since $(x,w)\in\E_2$, we know that $xu$ 
is not contained in the hyperplane $\E_{\leq 1}(w)$ and so $u$ is the only point on $xu$ contained in $\E_{\leq 1}(w)$. 

Symmetrically, for every point $z$ on $wy$ there is a point $z'$ on $xu$ such that $(z,z')\in\E_2$. Thus, $xu$ and $wy$ 
are opposite.\qed

\begin{cor} \label{maximal path}
For a point $x\in \E$, every $y\in\E_1(x)$ is collinear with a point $z\in\E_2(x)$.\qed
\end{cor}

\pf Let $u=[x,y]$. Applying Lemma \ref{opposite} (i) to $u$ and $y$ gives a point $z\in\E$ such that $y = [u,z]$. Applying 
Lemma \ref{opposite} (ii) to the path $xuyz$ gives that $(x,z)\in\E_2$.\qed

\begin{lem} \label{up and down}
Suppose $x\in\E$ and $\tau\in Y(x)$.
\begin{enumerate}
\item[\rm (i)] If $\tau$ fixes $y\in \E_1(x)$, then it also fixes every $z\in\E_2(x)$ that is collinear with $y$.
\item[\rm (ii)] If $\tau$ fixes $z\in\E_2(x)$, then it also fixes every $y\in\E_1(x)$ that is collinear with $z$.
\end{enumerate}
\end{lem}

\pf For the first claim, since $(x,y)\in\E_1$, the points $x$ and $y$ have a unique common neighbour $u=[x,y]$. Clearly, 
$(u,z)\in\E_1$ since $x\in\E_2(z)$. According to Lemma \ref{opposite} (i), there exists a point $v$ collinear with $z$ 
such that $(y,v)\in\E_1$. Note that since $\E_{\le 0}(y)$ is a subspace, every point in $zv\setminus\{z\}$ is in 
$\E_1(y)$. Hence we may choose $v$ to be the only point on $zv$ contained in $\E_1(x)$. Let $w$ the unique point collinear 
to both $x$ and $v$, that is, $w=[x,v]$.

Applying Lemma \ref{opposite} (ii) to the path $uyzv$ in the collinearity graph of $\E$ gives $(u,v)\in\E_2$. Note that 
$(w,z)\in\E_1$ (again, since $x\in\E_2(z)$) and $v=[w,z]$. Then applying Lemma \ref{opposite} (ii) to the path $yzvw$ 
gives $(y,w)\in\E_2$. So any two points at distance 3 (respectively, 2) in the hexagon $xuyzvw$ are in relation $\E_2$ 
(respectively, $\E_1$). 

Note that $vw$ is a cutting line with respect to $x$ and hence it is stabilised by $\tau$. Since $\tau$ fixes $y$ and 
stabilises $vw$, it fixes the unique point on $vw$ that is in the hyperplane $\E_{\leq 1}(y)$, namely, $\tau$ fixes $v$. 
Finally, $\tau$ fixes $z$ as it is the unique point collinear with both $y$ and $v$. This proves part (i).

For part (ii), we have that $u=[x,y]\in\E_{-1}(x)$ and so $\tau$ fixes $u$. As $u\in\E_1(z)$ and $y=[u,z]$, we also 
have that $\tau$ fixes $y$, as claimed.\qed

\medskip
We see that once $\tau$ fixes a point in $\E_1(x)\cup\E_2(x)$, it fixes many of its neighbours in this set. In order to 
fully exploit this, let us establish the following connectivity result.

\begin{lem} \label{connectivity}
Suppose $x\in\E$ and let $\Gamma$ be the graph on $\E_2(x)$ where two points are adjacent whenever they are at distance 
$1$ or $2$ in the collinearity graph of $\E$. Then $\Gamma$ is connected.
\end{lem}

\pf Suppose $y,z\in\E_2(x)$. If they are adjacent in $\Gamma$, then there is nothing to prove. Hence suppose that 
$(y,z)\in\E_2$. Let $yuvz$ be a shortest path from $y$ to $z$ in the collinearity graph of $\E$. If $u$ or $v$ is in 
$\E_2(x)$, then clearly $y$ and $z$ are in the same connected component of $\Gamma$. Hence we can assume that both $u$ and 
$v$ are in $\E_{\leq 1}(x)$.

By Definition \ref{def:rfs} (D), $(y,v)$ and $(u,z)$ are in $\E_1$. By Lemma \ref{opposite} (i), there exists $w$ 
collinear with $y$ such that $(u,w)\in\E_1$. By Lemma \ref{opposite} (iii), the lines $vz$ and $yw$ are opposite.
Note that by our assumption $v$ is the only point on $vz$ that is not in $\E_2(x)$. Let $s$ be the unique point of $yw$ 
contained in $\E_{\leq 1}(x)$. Since lines are thick, we can choose $t$ on $yw$ such that $t$ is not $y$ or $s$. Choose 
the unique point $t'$ on $vz$ such that $(t,t')\in\E_1$. This point exists since being in relation $\E_1$ is a bijection 
between the opposite lines $yw$ and $vz$. In particular, $t'$ is not $v$ and so $t'\in\E_2(x)$.  Now, $ytt'z$ is a path 
in $\Gamma$ between $y$ and $z$ and so $\Gamma$ is connected.\qed

\medskip
We now turn to the main property of point groups.

\begin{prop} \label{fixes}
For $x\in\E$ and $\tau\in Y(x)$, if $\tau$ fixes a point of $\E_1(x)$ then $\tau$ is the trivial automorphism.
That is, $Y(x)$ acts semiregularly on $\E_1(x)$.
%every cutting line with respect to $x$
\end{prop}

\pf Suppose $\tau$ fixes $y\in\E_1(x)$. By Corollary \ref{maximal path}, $y$ is collinear with a point $z\in\E_2(x)$. 
Lemma \ref{up and down} (i) implies that $\tau$ fixes every such $z$. Let $\Gamma$ be the graph on $\E_2(x)$ as 
defined in Lemma \ref{connectivity}. Suppose $z'$ is a neighbour of $z$ in $\Gamma$. If $z$ and $z'$ are 
collinear in $\E$, then they are both collinear with the unique point $y'\in\E_1(x)$. By Lemma \ref{up and down} (ii), 
$\tau$ fixes $y'$ and then also $\tau$ fixes $z'$ by Lemma \ref{up and down} (i). Similarly, if $z$ and $z'$ are not 
collinear in $\E$, then they must be at distance two in $\E$ and hence there exists a point $u$ collinear with both $z$ 
and $z'$. Whether $u\in\E_2(x)$ or $u\in\E_1(x)$, we see, as above, that $\tau$ fixes $u$ and then it fixes $z'$ as well. 
Hence whenever $\tau$ fixes $z\in\E_2(x)$, it also fixes all its neighbours in $\Gamma$. It now follows from the
connectivity property of $\Gamma$ given in Lemma  \ref{connectivity} that $\tau$ fixes every point in $\E_2(x)$. 

Finally, as every $y'\in\E_1(x)$ is collinear with some $z'\in\E_2(x)$, the automorphism $\tau$ fixes every such $y'$
by Lemma \ref{opposite} (ii). As $\tau\in Y(x)$, it also fixes every point in $\E_{\le 0}(x)$, and so we conclude 
that $\tau$ is trivial.\qed  

\medskip

Let us now assume that the root filtration space $\E$ is the extremal geometry of a Lie algebra $L$. Let $p$ be a point 
of $\E$, that is, $p$ is a 1-dimensional subspace of $L$ spanned by an extremal element $x$. According to 
\cite{ci1}, for every $\lambda\in F$, the map $\exp(x,\lambda)$ defined by $y\mapto y+\lambda[x,y]+\lambda^2g_x(y)x$ 
is an automorphism of $L$. Let $\Exp(p) = \{\exp(x,\lambda) \mid \lambda\in F\}$. 
It is easy to see that  $\exp(x,\lambda)\exp(x,\mu)=\exp(x,\lambda+\mu)$ for all $\lambda,\mu\in F$, which 
means that $\Exp(p)$ is a subgroup of $\aut(L)$. Moreover, as $\exp(\mu x,\lambda)=\exp(x,\mu\lambda)$, 
the definition of $\Exp(p)$ is independent of the choice of $x$ in $p$.

%Now assume that $L$ is a finite dimensional
%Lie algebra generated by its extremal elements with extremal geometry 
%$\E=\E(L)$ as in the hypothesis of our main result Theorem \ref{mainthm}.
%It is  our goal to show that the group $Exp(p)$, for $p\in\E$,  
%equals the group induced on $\E$ by a root subgroup of the 
%automorphism group of the corresponding building $\Delta$.
%Clearly, $Exp(p)$ is contained in such an induced group.

Consider the action of the elements of $\Exp(p)$ on the extremal geometry $\E$ of $L$.

\begin{lem} \label{mainLie} Let $p\in\E(L)$. Then
\begin{itemize}
\item[{\rm{(i)}}] $\Exp(p)$ acts trivially on $\E_{\le 0}(p)$, and
\item[{\rm{(ii)}}] $\Exp(p)$ acts transitively on $l\cap\E_1(p)$ for each cutting line $l$ with respect to $p$. 
\end{itemize}
\end{lem}

\pf If $q\in\E_{\le 0}(p)$, then $[x,y]=0$, where $x$ and $y$ are extremal elements spanning $p$ and $q$, 
respectively. By definition, $\exp(x,\lambda)$ fixes $y$, and so part (i) holds. 

Now suppose $l$ is a cutting line with respect to $p$ and $r\in\E_1(p)$ be a point on $l$. Let $q=[p,r]$. 
Select extremal elements $x$ and $z$ spanning $p$ and $r$, and set $y=[x,z]$. Then $y$ spans $q$. Furthermore, 
$l=qr=Fy+Fz$. The set $l\cap\E_1(p)$ equals $\{F(\lambda y+z)\mid\lambda\in F\}$. Since $r,q\in\E_{\le 1}(p)$, 
$g_x(z)=0=g_x(y)$ and, furthermore, $[x,y]=0$, as $q\in\E_{-1}(p)$. Therefore, $\exp(x,\mu)$ takes $\lambda y+z$ to 
$(\lambda y+z)+\mu[x,\lambda y+z]=\lambda y+z+\mu\lambda[x,y]+\mu[x,z]=(\mu+\lambda)y+z$. Clearly, this means 
that $\Exp(p)$ acts on $l\cap\E_1(p)$ and this action is transitive.
\qed

\medskip

%Clearly, $\Exp(p)$, being a group of automorphisms of $L$, acts on $\E=\E(L)$. 

\begin{cor} \label{regular}
For $p\in\E=\E(L)$, the image of $\Exp(p)$ in $\aut(\E)$ coincides with the point group $Y(p)$. Furthermore, 
the latter acts regularly on $l\cap\E_1(p)$ for each cutting line $l$ with respect to $p$.
\end{cor}

\pf By Lemma \ref{mainLie}, $\Exp(p)$ maps into $Y(p)$. Taking any $\tau\in Y(p)$ and a point 
$r\in l\cap\E_1(p)$, let $r'$ be the image of $r$ under $\tau$. It follows from Lemma \ref{mainLie} (ii) that 
the image of $\Exp(p)$ contains an automorphism $\tau'$ taking $r$ to the same point $r'$. Hence 
$\tau(\tau')^{-1}$ fixes $r$. However, by Proposition \ref{fixes}, $Y(p)$ is semiregular on $l\cap\E_1(p)$, 
which means that $\tau(\tau')^{-1}$ is the trivial automorphism, proving that $\tau'=\tau$. Therefore, 
the image of $\Exp(p)$ coincides with $Y(p)$. Also, since $Y(p)$ acts on $l\cap\E_1(p)$ semiregularly and 
transitively, we conclude that this action is regular.\qed

\medskip
Let us now switch to the language of buildings. Namely, we assume that the root filtration space $\E$ is 
the root shadow space of a thick irreducible spherical building $\Delta$ of simply-laced type with index set $I$. 
Recall from Section 
\ref{sec:rss} that the points of $\E$ are the $J$-shadows of $\Delta$ (where $J$ is the root set) and two 
points are collinear when the two corresponding $J$-shadows intersect a common $j$-panel for some $j\in J$. 
Let $\aut(\Delta)$ and $\aut^\circ(\Delta)$ denote the full groups of automorphisms and type-preserving automorphisms, 
respectively, of $\Delta$. Clearly, $\aut^\circ(\Delta)$ acts on the point set of $\E$ and, furthermore, it 
preserves the collinearity relation $\E_{-1}$. Thus $\aut^\circ(\Delta)$ induces a group of automorphisms of $\E$.

For a half-apartment $\alpha$ of $\Delta$, the \emph{root group} $U_\alpha$ is the subgroup of 
$\aut(\Delta)$ consisting of all those automorphisms that act trivially on any panel that intersects $\alpha$ 
in at least two chambers. Consequently, $U_\alpha$ acts trivially on $\alpha$. Selecting a chamber in the 
arctic region of $\alpha$, we see that any panel containing this chamber intersects $\alpha$ in two chambers, and so 
$U_\alpha$ stabilises a panel of each type in $I$. Hence $U_\alpha$ is type-preserving, that is, 
$U_\alpha\leq\aut^\circ(\Delta)$. 
In turn, this means that $U_\alpha$ stabilises setwise every residue that intersects $\alpha$. 

The building $\Delta$ is \emph{Moufang} if for each 
half-apartment $\alpha$ of $\Delta$ and for each panel $P$ containing just one chamber of  $\alpha$, the 
root group $U_\alpha$ acts transitively on $P\setminus P\cap\alpha$. It is well-known that all irreducible spherical 
buildings of rank at least $3$ are Moufang buildings (\cf Theorem 11.6 of \cite{w2}). Hence, in the simply-laced case, 
$\Delta$ is always Moufang unless it has type $A_2$. 

Suppose now that $\Delta$ is Moufang. Since $\Delta$ is of simply-laced type,
Theorem 30.14 of \cite{w1} implies that the 
root groups of $\Delta$ are all abelian. Combining this property with 
Proposition 3.16 in \cite{mw} we obtain the following result.

\begin{prop} \label{samePR=sameRG}
Suppose $\Delta$ is a Moufang building. If $\alpha$ and $\alpha'$ are two half-apartments with the same polar region, 
then $U_\alpha=U_{\alpha'}$.\qed
\end{prop}

This proposition means that, in the Moufang case, we can write $U_R$ instead of $U_\alpha$, where $R$ is the polar 
region of the half-apartment $\alpha$. For the next proposition, recall that we assume that $\E=\rsh(\Delta)$ 
and so the points of $\E$ are the $J$-shadows or, equivalently, polar regions of half-apartments.

\begin{prop} \label{mainBuilding} 
If $\Delta$ is Moufang, then for each point $R$ of $\E$, the image of $U_R$ in $\aut(\E)$ is 
contained in the point group $Y(R)$. 
\end{prop}

\pf Let $Q\in \E_{\le 0}(R)$ and choose an apartment $\Sigma$ that intersects both $R$ and $Q$. 
Let $\alpha$ and $\beta$ be the half-apartments in $\Sigma$, whose 
arctic regions are $R'=R\cap\Sigma$ and $Q'=Q\cap\Sigma$, respectively. 
Note that, by Proposition \ref{samePR=sameRG}, we have that $U_R=U_\alpha$. 
By Table \ref{tab:dict}, $\alpha$ has a nontrivial intersection with $Q'$, and hence with $Q$. Since $U_R=U_\alpha$ 
is type-preserving and it fixes a chamber from the $J$-shadow $Q$, it must stabilise $Q$ setwise and thus 
fix it as a point of $\E$.

Next, let $l$ be a cutting line with respect to $R$ and $Q\in\E_1(R)\cap l$. Again, select an apartment 
$\Sigma$ that intersects both $R$ and $Q$. Let $R'$, $Q'$, $\alpha$, and $\beta$ have the same meaning 
as in the paragraph above.
By Table \ref{tab:dict}, there is a unique half-apartment in $\Sigma$ whose arctic region $S'$ is collinear 
with both $R'$ and $Q'$ in $\rsh(\Sigma)$. Let $S$ be the corresponding polar region. Then $S=[R,Q]$. Since 
$S'$ and $Q'$ are collinear, they intersect a common $j$-panel of $\Sigma$, for some $j\in J$, that extends to
a unique $j$-panel $r$ of $\Delta$. The $j$-panel $r$ defines a line of $\E$ containing both $S$ and $Q$, which
coincides with the line $l$ since $\E$ is a partial linear space. Finally, 
since $U_R=U_\alpha$ is type-preserving and it fixes a chamber in $r$, it must stabilise $r$ and thus also $l$. 

Therefore the image of $U_R$ in $\aut(\E)$ lies in the point group $Y(R)$.\qed

\medskip
We previously mentioned that the only case where $\Delta$ can possibly be non-Moufang is when it has type 
$A_2$. The final proposition of this section shows that in our situation, where $\E$ arises as the 
extremal geometry of a Lie algebra, $\Delta$ is Moufang even for the type $A_2$.

\begin{prop} \label{A_2 is Moufang}
Suppose that the extremal geometry $\E(L)$ is isomorphic to the root shadow space of a building $\Delta$ 
of type $A_2$. Then $\Delta$ is a Moufang building.
\end{prop}

\pf Let $\E=\rsh(\Delta)$. Note that, for the type $A_2$, the root set $J$ and the index set $I$ coincide, 
namely, $I=J=\{1,2\}$. In particular, the points of $\E$ are the chambers and the lines are the panels of 
$\Delta$. Let $\alpha$ be a half-apartment of $\Delta$ whose polar region is the point (chamber) $p$ of 
$\E$. Then $\alpha$ contains three chambers including $p$.

Let $l$ be a panel containing two chambers of $\alpha$. Then one chamber must be $p$ and so $l\le\E_{-1}(p)$.
Hence the point group $Y(p)$ acts trivially on $l$. Secondly, let $l$ be a panel that contains exactly one 
chamber $q$ of $\alpha$. Then $q\ne p$. Since $\E_0$ is empty in the building of type $A_2$ (no two roots 
are perpendicular), we have that $l\setminus\{q\}\subset\E_1(p)$ and hence $l$ is a cutting 
line with respect to $p$. As $\E$ and $\E(L)$ are isomorphic, it follows from Corollary \ref{regular} that $Y(p)$ acts 
regularly (and hence transitively) on $l\setminus\{q\}$. This shows both that $Y(p)$ is contained in the root 
group $U_\alpha$ and that $\Delta$ is hence Moufang.\qed

%--------------------------------------- End: The Roots Subgroups --------------------------------------------------%

%--------------------------------------- Proof --------------------------------------------------%

\section{Proof of Theorem \ref{mainthm}}
\label{sec:proof}

In this section we complete the proof of the main result of the paper. We begin with a couple of technical
results about Moufang buildings.

Let $\Delta$ be an irreducible spherical building and assume it to be Moufang. Fix an apartment $\Sigma$ of 
$\Delta$. Let $U$ be the group generated by the root subgroups $\{ U_\alpha \mid \alpha\in\Sigma\}$.

\begin{lem} \label{apt^U=all} The building $\Delta$ is the union of $\{ \Sigma^g \mid g\in U\}$.
\end{lem}

\pf Let $c$ be a chamber in $\Delta$. We show that there exists $g\in U$ such that 
$c\in\Sigma^g$. We proceed by induction on $\dist(c, \Sigma) = \min_{x\in\Sigma} \dist(c,x)$ 
where $\dist(c,x)$ is the length of a minimal path between $c$ and $x$. 
Obviously if $\dist(x,\Sigma) = 0$, then choosing $g$ to be the identity suffices. 
Suppose now that $(x_0 = x, x_1, x_2, \ldots, x_{k-1}, x_k = c)$
is a minimal path of length $k$ from $\Sigma$ to $c$. In particular, $x\in\Sigma$ and 
$x_1\not\in\Sigma$. Let $P$ the panel containing $x$ and $x_1$. Let $x'$ be the unique 
chamber in $\Sigma$ different from $x$ that is contained in $P$. Let $\alpha$ be the unique 
half-apartment in $\Sigma$ containing $x$ but not $x'$. Since $\Delta$ is a Moufang building, there 
exists a $g\in U_\alpha$ such that $x_1^g = x'\in\Sigma$. Then 
$(x_1^g, x_2^g, \ldots, x_k^g = c^g)$ is a path of length less than $k$. By the inductive 
hypothesis, there exists $h\in U$ such that $c^g \in \Sigma^h$. Then $c\in\Sigma^{hg^{-1}}$.
But $hg^{-1}$ is in $U$ since $g\in U_\alpha$.\qed

\medskip

Recall that we identify $\rsh(\Sigma)$ with the corresponding set of points of $\rsh(\Delta)$ by identifying every 
arctic region in $\Sigma$ with the corresponding polar region in $\Delta$. Hence we can write under this identification 
that $\rsh(\Sigma) = \{ R\in\rsh(\Delta) \mid R\cap\Sigma\ne\emptyset \}$.

\begin{cor} \label{apt^U=all2} The root shadow space  $\rsh(\Delta)$ is the union of $\{{\rsh(\Sigma)}^g \mid g\in U\}$.
\end{cor}

\pf Let $R$ be in $\rsh(\Delta)$. Let  $c$ be a chamber of $R$. By Lemma \ref{apt^U=all} there exists $g\in U$ such 
that $c\in\Sigma^g$. 
But then $R\cap\Sigma^g$ is a nonempty intersection and thus $R\in{\rsh}(\Sigma^g) = {\rsh(\Sigma)}^g$,
as required.\qed

\medskip

We now embark upon the proof of Theorem \ref{mainthm}.
Let $L$ be a finite-dimensional simple Lie algebra that is generated by its extremal elements (or, equivalently, extremal points).
Let $\E=\E(L)$ be the extremal geometry of $L$ and assume it has lines. Then $\E$ is isomorphic to 
the root shadow space of an irreducible spherical building $\Delta$. We are additionally assuming 
that $\Delta$ is of simply laced type and so in view of Proposition \ref{A_2 is Moufang}, $\Delta$ is Moufang.
Let $L'$ be 
the subalgebra as defined in Section \ref{sec:subalg} with respect to a fixed apartment $\Sigma$ of 
$\Delta$. Namely, $L'$ is generated by the image $\E_\Sigma := \{ Fx_\alpha \mid \alpha\in\Phi\}$ of
$\rsh(\Sigma)$ under the isomorphism from $\rsh(\Delta)$ to $\E(L)$. 
Let $G$ be the group generated by
$\{ \Exp(x_\alpha) \mid \alpha\in\Phi \}$.
We state and prove the following observation. 

\begin{prop} \label{G is in aut(L')} The group $G$ leaves $L'$ invariant.
\end{prop}

\pf Select a root $\alpha$ in $\Phi$ and let $x\in L'$. Then 
\[ \exp(x_\alpha, t)x = x + t[x_\alpha,x]+ t^2 g_{x_\alpha}(x) x_\alpha \] 
and this is clearly in $L'$ by its construction. \qed

\medskip
By Corollary \ref{regular}, $G$ induces on $\E(L)$ the group generated by all the point groups 
$Y(Fx_\alpha)$. On the other hand, combining Proposition \ref{mainBuilding} with Corollary \ref{apt^U=all2}, we see that 
\[ \E = \bigcup_{g\in G} \E_\Sigma^g. \]
That is, $\E(L) = \bigcup_{g\in G} \{ Fx_\alpha \mid \alpha\in\Phi\}^g$. In particular, for each extremal point
$Fx$ in $\E(L)$, there exists an automorphism $g\in G$ such that $Fx = (Fx_\alpha)^g$.
But $(Fx_\alpha)^g$ is subspace of the subalgebra $L'$ by Proposition \ref{G is in aut(L')}. 
In particular, $L'$ contains all the extremal points from $\E(L)$. However, since $L$ is generated by 
its extremal points, we conclude that $L = L'$ and the proof of Theorem \ref{mainthm} is complete.

%--------------------------------------- End: Moufang Condition --------------------------------------------------%

%------------- Bibliography ---------------%

\end{document}